\documentstyle[amsfonts,12pt,a4]{article}
\makeatletter
\newcount\@hour\newcount\@minute\chardef\@x10\chardef\@xv60
\def\tcitime{
\def\@time{%
  \@minute\time\@hour\@minute\divide\@hour\@xv
  \ifnum\@hour<\@x 0\fi\the\@hour:%
  \multiply\@hour\@xv\advance\@minute-\@hour
  \ifnum\@minute<\@x 0\fi\the\@minute
  }}%
\def\QCTOpt[#1]#2{%
  \def\QCTOptB{#1}
  \def\QCTOptA{#2}
}
\def\QCTNOpt#1{%
  \def\QCTOptA{#1}
  \let\QCTOptB\empty
}
\def\Qct{%
  \@ifnextchar[{%
    \QCTOpt}{\QCTNOpt}
}
\def\QCBOpt[#1]#2{%
  \def\QCBOptB{#1}
  \def\QCBOptA{#2}
}
\def\QCBNOpt#1{%
  \def\QCBOptA{#1}
  \let\QCBOptB\empty
}
\def\Qcb{%
  \@ifnextchar[{%
    \QCBOpt}{\QCBNOpt}
}
\def\PrepCapArgs{%
  \ifx\QCBOptA\empty
    \ifx\QCTOptA\empty
      {}%
    \else
      \ifx\QCTOptB\empty
        {\QCTOptA}%
      \else
        [\QCTOptB]{\QCTOptA}%
      \fi
    \fi
  \else
    \ifx\QCBOptA\empty
      {}%
    \else
      \ifx\QCBOptB\empty
        {\QCBOptA}%
      \else
        [\QCBOptB]{\QCBOptA}%
      \fi
    \fi
  \fi
}
\newcount\GRAPHICSTYPE
\GRAPHICSTYPE=\z@
\def\GRAPHICSPS#1{%
 \ifcase\GRAPHICSTYPE
   \special{ps: #1}%
 \or
   \special{language "PS", include "#1"}%
 \fi
}%
\def\graffile#1#2#3#4{%
    \leavevmode
    \raise -#4 \BOXTHEFRAME{%
        \hbox to #2{\raise #3\hbox{\null #1}}}%
}%
\def\draftbox#1#2#3#4{%
 \leavevmode\raise -#4 \hbox{%
  \frame{\rlap{\protect\tiny #1}\hbox to #2%
   {\vrule height#3 width\z@ depth\z@\hfil}%
  }%
 }%
}%
\newcount\draft
\draft=\z@

\newif\ifwasdraft
\wasdraftfalse

\def\GRAPHIC#1#2#3#4#5{%
 \ifnum\draft=\@ne\draftbox{#2}{#3}{#4}{#5}%
  \else\graffile{#1}{#3}{#4}{#5}%
  \fi
 }%
\def\addtoLaTeXparams#1{%
    \edef\LaTeXparams{\LaTeXparams #1}}%
\newif\ifBoxFrame \BoxFramefalse
\newif\ifOverFrame \OverFramefalse
\newif\ifUnderFrame \UnderFramefalse

\def\BOXTHEFRAME#1{%
   \hbox{%
      \ifBoxFrame
         \frame{#1}%
      \else
         {#1}%
      \fi
   }%
}

\def\doFRAMEparams#1{\BoxFramefalse\OverFramefalse\UnderFramefalse\readFRAMEparams#1\end}%
\def\readFRAMEparams#1{%
 \ifx#1\end%
  \let\next=\relax
  \else
  \ifx#1i\dispkind=\z@\fi
  \ifx#1d\dispkind=\@ne\fi
  \ifx#1f\dispkind=\tw@\fi
  \ifx#1t\addtoLaTeXparams{t}\fi
  \ifx#1b\addtoLaTeXparams{b}\fi
  \ifx#1p\addtoLaTeXparams{p}\fi
  \ifx#1h\addtoLaTeXparams{h}\fi
  \ifx#1X\BoxFrametrue\fi
  \ifx#1O\OverFrametrue\fi
  \ifx#1U\UnderFrametrue\fi
  \ifx#1w
    \ifnum\draft=1\wasdrafttrue\else\wasdraftfalse\fi
    \draft=\@ne
  \fi
  \let\next=\readFRAMEparams
  \fi
 \next
 }%
\def\IFRAME#1#2#3#4#5#6{%
      \bgroup
      \let\QCTOptA\empty
      \let\QCTOptB\empty
      \let\QCBOptA\empty
      \let\QCBOptB\empty
      #6%
      \parindent=0pt%
      \leftskip=0pt
      \rightskip=0pt
      \setbox0 = \hbox{\QCBOptA}%
      \@tempdima = #1\relax
      \ifOverFrame
          \typeout{This is not implemented yet}%
          \show\HELP
      \else
         \ifdim\wd0>\@tempdima
            \advance\@tempdima by \@tempdima
            \ifdim\wd0 >\@tempdima
               \textwidth=\@tempdima
               \setbox1 =\vbox{%
                  \noindent\hbox to \@tempdima{\hfill\GRAPHIC{#5}{#4}{#1}{#2}{#3}\hfill}\\%
                  \noindent\hbox to \@tempdima{\parbox[b]{\@tempdima}{\QCBOptA}}%
               }%
               \wd1=\@tempdima
            \else
               \textwidth=\wd0
               \setbox1 =\vbox{%
                 \noindent\hbox to \wd0{\hfill\GRAPHIC{#5}{#4}{#1}{#2}{#3}\hfill}\\%
                 \noindent\hbox{\QCBOptA}%
               }%
               \wd1=\wd0
            \fi
         \else
            \ifdim\wd0>0pt
              \hsize=\@tempdima
              \setbox1 =\vbox{%
                \unskip\GRAPHIC{#5}{#4}{#1}{#2}{0pt}%
                \break
                \unskip\hbox to \@tempdima{\hfill \QCBOptA\hfill}%
              }%
              \wd1=\@tempdima
           \else
              \hsize=\@tempdima
              \setbox1 =\vbox{%
                \unskip\GRAPHIC{#5}{#4}{#1}{#2}{0pt}%
              }%
              \wd1=\@tempdima
           \fi
         \fi
         \@tempdimb=\ht1
         \advance\@tempdimb by \dp1
         \advance\@tempdimb by -#2%
         \advance\@tempdimb by #3%
         \leavevmode
         \raise -\@tempdimb \hbox{\box1}%
      \fi
      \egroup%
}%
\def\DFRAME#1#2#3#4#5{%
 \begin{center}
     \let\QCTOptA\empty
     \let\QCTOptB\empty
     \let\QCBOptA\empty
     \let\QCBOptB\empty
     \ifOverFrame 
        #5\QCTOptA\par
     \fi
     \GRAPHIC{#4}{#3}{#1}{#2}{\z@}
     \ifUnderFrame 
        \par #5\QCBOptA
     \fi
 \end{center}%
 }%
\def\FFRAME#1#2#3#4#5#6#7{%
 \begin{figure}[#1]%
  \let\QCTOptA\empty
  \let\QCTOptB\empty
  \let\QCBOptA\empty
  \let\QCBOptB\empty
  \ifOverFrame
    #4
    \ifx\QCTOptA\empty
    \else
      \ifx\QCTOptB\empty
        \caption{\QCTOptA}%
      \else
        \caption[\QCTOptB]{\QCTOptA}%
      \fi
    \fi
    \ifUnderFrame\else
      \label{#5}%
    \fi
  \else
    \UnderFrametrue%
  \fi
  \begin{center}\GRAPHIC{#7}{#6}{#2}{#3}{\z@}\end{center}%
  \ifUnderFrame
    #4
    \ifx\QCBOptA\empty
      \caption{}%
    \else
      \ifx\QCBOptB\empty
        \caption{\QCBOptA}%
      \else
        \caption[\QCBOptB]{\QCBOptA}%
      \fi
    \fi
    \label{#5}%
  \fi
  \end{figure}%
 }%
\newcount\dispkind%
\def\FRAME#1#2#3#4#5#6#7#8{%
 \ifnum\draft=\@ne
   \wasdrafttrue
 \else
   \wasdraftfalse%
 \fi
 \def\LaTeXparams{}%
 \dispkind=\z@
 \def\LaTeXparams{}%
 \doFRAMEparams{#1}%
 \ifnum\dispkind=\z@\IFRAME{#2}{#3}{#4}{#7}{#8}{#5}\else
  \ifnum\dispkind=\@ne\DFRAME{#2}{#3}{#7}{#8}{#5}\else
   \ifnum\dispkind=\tw@
    \edef\@tempa{\noexpand\FFRAME{\LaTeXparams}}%
    \@tempa{#2}{#3}{#5}{#6}{#7}{#8}%
    \fi
   \fi
  \fi
  \ifwasdraft\draft=1\else\draft=0\fi{}%
 }%
\def\TEXUX#1{"texux"}

\long\def\QQQ#1#2{%
     \long\expandafter\def\csname#1\endcsname{#2}}%
\@ifundefined{QTP}{\def\QTP#1{}}{}
\@ifundefined{Qlb}{}{}
\@ifundefined{Qlt}{}{}
\long\def\QQA#1#2{}%
\def\QTR#1#2{{\csname#1\endcsname #2}}
\def\EXPAND#1[#2]#3{}%
\def\NOEXPAND#1[#2]#3{}%
\def\LaTeXparent#1{}%
\def\ChildStyles#1{}%
\def\ChildDefaults#1{}%
\def\QTagDef#1#2#3{}%
\def\QQfnmark#1{\footnotemark}

\@ifundefined{INDEX}{\def\INDEX#1#2{}{}}{}%
\@ifundefined{SUBINDEX}{\def\SUBINDEX#1#2#3{}{}{}}{}%
\def\initial#1{\bigbreak{\raggedright\large\bf #1}\kern 2\p@
   \penalty3000}%
\@ifundefined{ZZZ}{}{\makeatletter\input gnuindex.sty\makeatother\makeindex\makeatletter}%
\@ifundefined{abstract}{%
 \def\abstract{%
  \if@twocolumn
   \section*{Abstract (Not appropriate in this style!)}%
   \else \small 
   \begin{center}{\bf Abstract\vspace{-.5em}\vspace{\z@}}\end{center}%
   \quotation 
   \fi
  }%
 }{%
 }%
\@ifundefined{endabstract}{\def\endabstract
  {\if@twocolumn\else\endquotation\fi}}{}%
\@ifundefined{maketitle}{\def\maketitle#1{}}{}%
\@ifundefined{affiliation}{\def\affiliation#1{}}{}%
\@ifundefined{proof}{}{}%
\@ifundefined{endproof}{}{}%
\@ifundefined{newfield}{\def\newfield#1#2{}}{}%
\@ifundefined{chapter}{\def\chapter#1{\par(Chapter head:)#1\par }%
 \newcount\c@chapter}{}%
\@ifundefined{part}{\def\part#1{\par(Part head:)#1\par }}{}%
\@ifundefined{section}{\def\section#1{\par(Section head:)#1\par }}{}%
\@ifundefined{subsection}{\def\subsection#1%
 {\par(Subsection head:)#1\par }}{}%
\@ifundefined{subsubsection}{\def\subsubsection#1%
 {\par(Subsubsection head:)#1\par }}{}%
\@ifundefined{paragraph}{\def\paragraph#1%
 {\par(Subsubsubsection head:)#1\par }}{}%
\@ifundefined{subparagraph}{\def\subparagraph#1%
 {\par(Subsubsubsubsection head:)#1\par }}{}%
\@ifundefined{therefore}{}{}%
\@ifundefined{backepsilon}{}{}%
\@ifundefined{yen}{}{}%
\@ifundefined{registered}{%
   \def\registered{\relax\ifmmode{}\r@gistered
                    \else$\m@th\r@gistered$\fi}%
 \def\r@gistered{^{\ooalign
  {\hfil\raise.07ex\hbox{$\scriptstyle\rm\text{R}$}\hfil\crcr
  \mathhexbox20D}}}}{}%
\@ifundefined{Eth}{}{}%
\@ifundefined{eth}{}{}%
\@ifundefined{Thorn}{}{}%
\@ifundefined{thorn}{}{}%
\@ifundefined{degree}{}{}%
\newdimen\theight
\def\Column{%
 \vadjust{\setbox\z@=\hbox{\scriptsize\quad\quad tcol}%
  \theight=\ht\z@\advance\theight by \dp\z@\advance\theight by \lineskip
  \kern -\theight \vbox to \theight{%
   \rightline{\rlap{\box\z@}}%
   \vss
   }%
  }%
 }%
\def\qed{%
 \ifhmode\unskip\nobreak\fi\ifmmode\ifinner\else\hskip5\p@\fi\fi
 \hbox{\hskip5\p@\vrule width4\p@ height6\p@ depth1.5\p@\hskip\p@}%
 }%
\def\miss{\hbox{\vrule height2\p@ width 2\p@ depth\z@}}%
\def\tcol#1{{\baselineskip=6\p@ \vcenter{#1}} \Column}  %
\def\newfmtname{LaTeX2e}
\def\chkcompat{%
   \if@compatibility
   \else
     \usepackage{latexsym}
   \fi
}

\ifx\fmtname\newfmtname
  \DeclareOldFontCommand{\rm}{\normalfont\rmfamily}{\mathrm}
  \DeclareOldFontCommand{\sf}{\normalfont\sffamily}{\mathsf}
  \DeclareOldFontCommand{\tt}{\normalfont\ttfamily}{\mathtt}
  \DeclareOldFontCommand{\bf}{\normalfont\bfseries}{\mathbf}
  \DeclareOldFontCommand{\it}{\normalfont\itshape}{\mathit}
  \DeclareOldFontCommand{\sl}{\normalfont\slshape}{\@nomath\sl}
  \DeclareOldFontCommand{\sc}{\normalfont\scshape}{\@nomath\sc}
  \chkcompat
\fi

\def\alpha{\Greekmath 010B }%
\def\beta{\Greekmath 010C }%
\def\gamma{\Greekmath 010D }%
\def\delta{\Greekmath 010E }%
\def\epsilon{\Greekmath 010F }%
\def\zeta{\Greekmath 0110 }%
\def\eta{\Greekmath 0111 }%
\def\theta{\Greekmath 0112 }%
\def\iota{\Greekmath 0113 }%
\def\kappa{\Greekmath 0114 }%
\def\lambda{\Greekmath 0115 }%
\def\mu{\Greekmath 0116 }%
\def\nu{\Greekmath 0117 }%
\def\xi{\Greekmath 0118 }%
\def\pi{\Greekmath 0119 }%
\def\rho{\Greekmath 011A }%
\def\sigma{\Greekmath 011B }%
\def\tau{\Greekmath 011C }%
\def\upsilon{\Greekmath 011D }%
\def\phi{\Greekmath 011E }%
\def\chi{\Greekmath 011F }%
\def\psi{\Greekmath 0120 }%
\def\omega{\Greekmath 0121 }%
\def\varepsilon{\Greekmath 0122 }%
\def\vartheta{\Greekmath 0123 }%
\def\varpi{\Greekmath 0124 }%
\def\varrho{\Greekmath 0125 }%
\def\varsigma{\Greekmath 0126 }%
\def\varphi{\Greekmath 0127 }%

\def\nabla{\Greekmath 0272 }

\def\Greekmath#1#2#3#4{%
    \if@compatibility
        \ifnum\mathgroup=\symbold
           \mathchoice{\mbox{\boldmath$\displaystyle\mathchar"#1#2#3#4$}}%
                      {\mbox{\boldmath$\textstyle\mathchar"#1#2#3#4$}}%
                      {\mbox{\boldmath$\scriptstyle\mathchar"#1#2#3#4$}}%
                      {\mbox{\boldmath$\scriptscriptstyle\mathchar"#1#2#3#4$}}%
        \else
           \mathchar"#1#2#3#4%
        \fi 
    \else 
        \ifnum\mathgroup=5 
           \mathchoice{\mbox{\boldmath$\displaystyle\mathchar"#1#2#3#4$}}%
                      {\mbox{\boldmath$\textstyle\mathchar"#1#2#3#4$}}%
                      {\mbox{\boldmath$\scriptstyle\mathchar"#1#2#3#4$}}%
                      {\mbox{\boldmath$\scriptscriptstyle\mathchar"#1#2#3#4$}}%
        \else
           \mathchar"#1#2#3#4%
        \fi     	    
	  \fi}

\newif\ifGreekBold  \GreekBoldfalse
\let\SAVEPBF=\pbf
\def\pbf{\GreekBoldtrue\SAVEPBF}%

\@ifundefined{theorem}{}{}
\@ifundefined{lemma}{}{}
\@ifundefined{corollary}{}{}
\@ifundefined{conjecture}{}{}
\@ifundefined{proposition}{}{}
\@ifundefined{axiom}{}{}
\@ifundefined{remark}{}{}
\@ifundefined{example}{}{}
\@ifundefined{exercise}{}{}
\@ifundefined{definition}{}{}

\@ifundefined{mathletters}{%
  \newcounter{equationnumber}  
  \def\mathletters{%
     \addtocounter{equation}{1}
     \edef\@currentlabel{\theequation}%
     \setcounter{equationnumber}{\c@equation}
     \setcounter{equation}{0}%
     \edef\theequation{\@currentlabel\noexpand\alph{equation}}%
  }
  
}{}

\@ifundefined{BibTeX}{%
    \def\BibTeX{{\rm B\kern-.05em{\sc i\kern-.025em b}\kern-.08em
                 T\kern-.1667em\lower.7ex\hbox{E}\kern-.125emX}}}{}%
\@ifundefined{AmS}%
    {\def\AmS{{\protect\usefont{OMS}{cmsy}{m}{n}%
                A\kern-.1667em\lower.5ex\hbox{M}\kern-.125emS}}}{}%
\@ifundefined{AmSTeX}{}{}%
\ifx\ds@amstex\relax
\makeatother\endinput
\else
   \@ifpackageloaded{amstex}%
      {\message{amstex already loaded}\makeatother\endinput}
      {}
\fi
\let\DOTSI\relax
\def\RIfM@{\relax\ifmmode}%
\def\FN@{\futurelet\next}%
\newcount\intno@
\def\iint{\DOTSI\intno@\tw@\FN@\ints@}%
\def\iiint{\DOTSI\intno@\thr@@\FN@\ints@}%
\def\iiiint{\DOTSI\intno@4 \FN@\ints@}%
\def\idotsint{\DOTSI\intno@\z@\FN@\ints@}%
\def\ints@{\findlimits@\ints@@}%
\newif\iflimtoken@
\newif\iflimits@
\def\findlimits@{\limtoken@true\ifx\next\limits\limits@true
 \else\ifx\next\nolimits\limits@false\else
 \limtoken@false\ifx\ilimits@\nolimits\limits@false\else
 \ifinner\limits@false\else\limits@true\fi\fi\fi\fi}%
\def\multint@{\int\ifnum\intno@=\z@\intdots@                          
 \else\intkern@\fi                                                    
 \ifnum\intno@>\tw@\int\intkern@\fi                                   
 \ifnum\intno@>\thr@@\int\intkern@\fi                                 
 \int}
\def\multintlimits@{\intop\ifnum\intno@=\z@\intdots@\else\intkern@\fi
 \ifnum\intno@>\tw@\intop\intkern@\fi
 \ifnum\intno@>\thr@@\intop\intkern@\fi\intop}%
\def\intic@{%
    \mathchoice{\hskip.5em}{\hskip.4em}{\hskip.4em}{\hskip.4em}}%
\def\negintic@{\mathchoice
 {\hskip-.5em}{\hskip-.4em}{\hskip-.4em}{\hskip-.4em}}%
\def\ints@@{\iflimtoken@                                              
 \def\ints@@@{\iflimits@\negintic@
   \mathop{\intic@\multintlimits@}\limits                             
  \else\multint@\nolimits\fi                                          
  \eat@}
 \else                                                                
 \def\ints@@@{\iflimits@\negintic@
  \mathop{\intic@\multintlimits@}\limits\else
  \multint@\nolimits\fi}\fi\ints@@@}%
\def\intkern@{\mathchoice{\!\!\!}{\!\!}{\!\!}{\!\!}}%
\def\plaincdots@{\mathinner{\cdotp\cdotp\cdotp}}%
\def\intdots@{\mathchoice{\plaincdots@}%
 {{\cdotp}\mkern1.5mu{\cdotp}\mkern1.5mu{\cdotp}}%
 {{\cdotp}\mkern1mu{\cdotp}\mkern1mu{\cdotp}}%
 {{\cdotp}\mkern1mu{\cdotp}\mkern1mu{\cdotp}}}%
\def\RIfM@{\relax\protect\ifmmode}
\def\text{\RIfM@\expandafter\text@\else\expandafter\mbox\fi}
\let\nfss@text\text
\def\text@#1{\mathchoice
   {\textdef@\displaystyle\f@size{#1}}%
   {\textdef@\textstyle\tf@size{\firstchoice@false #1}}%
   {\textdef@\textstyle\sf@size{\firstchoice@false #1}}%
   {\textdef@\textstyle \ssf@size{\firstchoice@false #1}}%
   \glb@settings}

\def\textdef@#1#2#3{\hbox{{%
                    \everymath{#1}%
                    \let\f@size#2\selectfont
                    #3}}}
\newif\iffirstchoice@
\firstchoice@true
\def\Let@{\relax\iffalse{\fi\let\\=\cr\iffalse}\fi}%
\def\vspace@{\def\vspace##1{\crcr\noalign{\vskip##1\relax}}}%
\def\multilimits@{\bgroup\vspace@\Let@
 \baselineskip\fontdimen10 \scriptfont\tw@
 \advance\baselineskip\fontdimen12 \scriptfont\tw@
 \lineskip\thr@@\fontdimen8 \scriptfont\thr@@
 \lineskiplimit\lineskip
 \vbox\bgroup\ialign\bgroup\hfil$\m@th\scriptstyle{##}$\hfil\crcr}%
\def\Sb{_\multilimits@}%
\def\endSb{\crcr\egroup\egroup\egroup}%
\def\Sp{^\multilimits@}%

\newdimen\ex@
\def\rightarrowfill@#1{$#1\m@th\mathord-\mkern-6mu\cleaders
 \hbox{$#1\mkern-2mu\mathord-\mkern-2mu$}\hfill
 \mkern-6mu\mathord\rightarrow$}%
\def\leftarrowfill@#1{$#1\m@th\mathord\leftarrow\mkern-6mu\cleaders
 \hbox{$#1\mkern-2mu\mathord-\mkern-2mu$}\hfill\mkern-6mu\mathord-$}%
\def\leftrightarrowfill@#1{$#1\m@th\mathord\leftarrow
\mkern-6mu\cleaders
 \hbox{$#1\mkern-2mu\mathord-\mkern-2mu$}\hfill
 \mkern-6mu\mathord\rightarrow$}%
\def\overrightarrow{\mathpalette\overrightarrow@}%
\def\overrightarrow@#1#2{\vbox{\ialign{##\crcr\rightarrowfill@#1\crcr
 \noalign{\kern-\ex@\nointerlineskip}$\m@th\hfil#1#2\hfil$\crcr}}}%

\def\overleftarrow{\mathpalette\overleftarrow@}%
\def\overleftarrow@#1#2{\vbox{\ialign{##\crcr\leftarrowfill@#1\crcr
 \noalign{\kern-\ex@\nointerlineskip}$\m@th\hfil#1#2\hfil$\crcr}}}%
\def\overleftrightarrow{\mathpalette\overleftrightarrow@}%
\def\overleftrightarrow@#1#2{\vbox{\ialign{##\crcr
   \leftrightarrowfill@#1\crcr
 \noalign{\kern-\ex@\nointerlineskip}$\m@th\hfil#1#2\hfil$\crcr}}}%
\def\underrightarrow{\mathpalette\underrightarrow@}%
\def\underrightarrow@#1#2{\vtop{\ialign{##\crcr$\m@th\hfil#1#2\hfil
  $\crcr\noalign{\nointerlineskip}\rightarrowfill@#1\crcr}}}%

\def\underleftarrow{\mathpalette\underleftarrow@}%
\def\underleftarrow@#1#2{\vtop{\ialign{##\crcr$\m@th\hfil#1#2\hfil
  $\crcr\noalign{\nointerlineskip}\leftarrowfill@#1\crcr}}}%
\def\underleftrightarrow{\mathpalette\underleftrightarrow@}%
\def\underleftrightarrow@#1#2{\vtop{\ialign{##\crcr$\m@th
  \hfil#1#2\hfil$\crcr
 \noalign{\nointerlineskip}\leftrightarrowfill@#1\crcr}}}%

\def\qopnamewl@#1{\mathop{\operator@font#1}\nlimits@}
\let\nlimits@\displaylimits
\def\setboxz@h{\setbox\z@\hbox}

\def\varlim@#1#2{\mathop{\vtop{\ialign{##\crcr
 \hfil$#1\m@th\operator@font lim$\hfil\crcr
 \noalign{\nointerlineskip}#2#1\crcr
 \noalign{\nointerlineskip\kern-\ex@}\crcr}}}}

 \def\rightarrowfill@#1{\m@th\setboxz@h{$#1-$}\ht\z@\z@
  $#1\copy\z@\mkern-6mu\cleaders
  \hbox{$#1\mkern-2mu\box\z@\mkern-2mu$}\hfill
  \mkern-6mu\mathord\rightarrow$}
\def\leftarrowfill@#1{\m@th\setboxz@h{$#1-$}\ht\z@\z@
  $#1\mathord\leftarrow\mkern-6mu\cleaders
  \hbox{$#1\mkern-2mu\copy\z@\mkern-2mu$}\hfill
  \mkern-6mu\box\z@$}

\def\projlim{\qopnamewl@{proj\,lim}}
\def\injlim{\qopnamewl@{inj\,lim}}
\def\varinjlim{\mathpalette\varlim@\rightarrowfill@}
\def\varprojlim{\mathpalette\varlim@\leftarrowfill@}
\def\varliminf{\mathpalette\varliminf@{}}
\def\varliminf@#1{\mathop{\underline{\vrule\@depth.2\ex@\@width\z@
   \hbox{$#1\m@th\operator@font lim$}}}}
\def\varlimsup{\mathpalette\varlimsup@{}}
\def\varlimsup@#1{\mathop{\overline
  {\hbox{$#1\m@th\operator@font lim$}}}}

\def\dfrac#1#2{{\displaystyle {#1 \over #2}}}%
\begingroup \catcode `|=0 \catcode `[= 1
\catcode`]=2 \catcode `\{=12 \catcode `\}=12
\catcode`\\=12 
|gdef|@alignverbatim#1\end{align}[#1|end[align]]
|gdef|@salignverbatim#1\end{align*}[#1|end[align*]]

|gdef|@alignatverbatim#1\end{alignat}[#1|end[alignat]]
|gdef|@salignatverbatim#1\end{alignat*}[#1|end[alignat*]]

|gdef|@xalignatverbatim#1\end{xalignat}[#1|end[xalignat]]
|gdef|@sxalignatverbatim#1\end{xalignat*}[#1|end[xalignat*]]

|gdef|@gatherverbatim#1\end{gather}[#1|end[gather]]
|gdef|@sgatherverbatim#1\end{gather*}[#1|end[gather*]]

|gdef|@gatherverbatim#1\end{gather}[#1|end[gather]]
|gdef|@sgatherverbatim#1\end{gather*}[#1|end[gather*]]

|gdef|@multilineverbatim#1\end{multiline}[#1|end[multiline]]
|gdef|@smultilineverbatim#1\end{multiline*}[#1|end[multiline*]]

|gdef|@arraxverbatim#1\end{arrax}[#1|end[arrax]]
|gdef|@sarraxverbatim#1\end{arrax*}[#1|end[arrax*]]

|gdef|@tabulaxverbatim#1\end{tabulax}[#1|end[tabulax]]
|gdef|@stabulaxverbatim#1\end{tabulax*}[#1|end[tabulax*]]

|endgroup

\def\align{\@verbatim \frenchspacing\@vobeyspaces \@alignverbatim
You are using the "align" environment in a style in which it is not defined.}

\@namedef{align*}{\@verbatim\@salignverbatim
You are using the "align*" environment in a style in which it is not defined.}
\expandafter\let\csname endalign*\endcsname =\endtrivlist

\def\alignat{\@verbatim \frenchspacing\@vobeyspaces \@alignatverbatim
You are using the "alignat" environment in a style in which it is not defined.}

\@namedef{alignat*}{\@verbatim\@salignatverbatim
You are using the "alignat*" environment in a style in which it is not defined.}
\expandafter\let\csname endalignat*\endcsname =\endtrivlist

\def\xalignat{\@verbatim \frenchspacing\@vobeyspaces \@xalignatverbatim
You are using the "xalignat" environment in a style in which it is not defined.}

\@namedef{xalignat*}{\@verbatim\@sxalignatverbatim
You are using the "xalignat*" environment in a style in which it is not defined.}
\expandafter\let\csname endxalignat*\endcsname =\endtrivlist

\def\gather{\@verbatim \frenchspacing\@vobeyspaces \@gatherverbatim
You are using the "gather" environment in a style in which it is not defined.}

\@namedef{gather*}{\@verbatim\@sgatherverbatim
You are using the "gather*" environment in a style in which it is not defined.}
\expandafter\let\csname endgather*\endcsname =\endtrivlist

\def\multiline{\@verbatim \frenchspacing\@vobeyspaces \@multilineverbatim
You are using the "multiline" environment in a style in which it is not defined.}

\@namedef{multiline*}{\@verbatim\@smultilineverbatim
You are using the "multiline*" environment in a style in which it is not defined.}
\expandafter\let\csname endmultiline*\endcsname =\endtrivlist

\def\arrax{\@verbatim \frenchspacing\@vobeyspaces \@arraxverbatim
You are using a type of "array" construct that is only allowed in AmS-LaTeX.}

\def\tabulax{\@verbatim \frenchspacing\@vobeyspaces \@tabulaxverbatim
You are using a type of "tabular" construct that is only allowed in AmS-LaTeX.}

\@namedef{arrax*}{\@verbatim\@sarraxverbatim
You are using a type of "array*" construct that is only allowed in AmS-LaTeX.}
\expandafter\let\csname endarrax*\endcsname =\endtrivlist

\@namedef{tabulax*}{\@verbatim\@stabulaxverbatim
You are using a type of "tabular*" construct that is only allowed in AmS-LaTeX.}
\expandafter\let\csname endtabulax*\endcsname =\endtrivlist

\def\@@eqncr{\let\@tempa\relax
    \ifcase\@eqcnt \def\@tempa{& & &}\or \def\@tempa{& &}%
      \else \def\@tempa{&}\fi
     \@tempa
     \if@eqnsw
        \iftag@
           \@taggnum
        \else
           \@eqnnum\stepcounter{equation}%
        \fi
     \fi
     \global\tag@false
     \global\@eqnswtrue
     \global\@eqcnt\z@\cr}

 \def\endequation{%
     \ifmmode\ifinner 
      \iftag@
        \addtocounter{equation}{-1} 
        $\hfil
           \displaywidth\linewidth\@taggnum\egroup \endtrivlist
        \global\tag@false
        \global\@ignoretrue   
      \else
        $\hfil
           \displaywidth\linewidth\@eqnnum\egroup \endtrivlist
        \global\tag@false
        \global\@ignoretrue 
      \fi
     \else   
      \iftag@
        \addtocounter{equation}{-1} 
        \eqno \hbox{\@taggnum}
        \global\tag@false%
        $$\global\@ignoretrue
      \else
        \eqno \hbox{\@eqnnum}
        $$\global\@ignoretrue
      \fi
     \fi\fi
 } 

 \newif\iftag@ \tag@false
 
 \def\tag{\@ifnextchar*{\@tagstar}{\@tag}}
 \def\@tag#1{%
     \global\tag@true
     \global\def\@taggnum{(#1)}}
 \def\@tagstar*#1{%
     \global\tag@true
     \global\def\@taggnum{#1}%
}

\makeatother

\QQQ{Language}{American English}
\def\ba{\begin{array}}
\def\ea{\end{array}}
\def\beq{\begin{equation}}
\def\eeq{\end{equation}}
\def\bea{\begin{eqnarray}}
\def\eea{\end{eqnarray}}
\def\si{\sigma}
\def\tsi{\tilde{\sigma}}
\def\vf{\varphi}
\def\ve{\varepsilon}
\def\om{\omega}
\def\Om{\Omega}
\def\Ga{\Gamma}
\def\lam{\lambda}
\def\Lam{\Lambda}
\def\De{{\rm Det}\; }
\def\Pf{{\rm Pfaff}\; }
\def\Tr{{\rm Tr}\; }

\begin{document}

\newcommand{\norm}[1]{{\protect\normalsize{#1}}}
\newcommand{\LAP}
{{\small E}\norm{N}{\large S}{\Large L}{\large A}\norm{P}{\small P}}
\renewcommand{\thefootnote}{\fnsymbol{footnote}} \newpage
\pagestyle{empty} \setcounter{page}{0}

\begin{center}
{\Large {\bf REPRESENTATIVES OF THE THOM CLASS\\[0.24cm]OF A VECTOR
BUNDLE$^{\dagger}$}}%
\\[1cm]

\vspace{10mm}

{\large Michel Bauer$^{1}$ and Frank Thuillier$^{2}$}

{\em Laboratoire de Physique Th\'eorique }{\small E}{{\normalsize {N}}}%
{\large S}{\Large L}{\large A}{{\normalsize {P}}}{\small P}\footnote{%
URA 14-36 du CNRS, associ\'ee \`a l'Ecole Normale Sup\'erieure de Lyon,
\par
\noindent
$^{\dagger}$ Work partially supported by European Community Contract
ERBCHRXCT920069
\par
\noindent 
$^1$ Groupe d'Annecy: LAPP, Chemin de Bellevue BP 110, F-74941
Annecy-le-Vieux Cedex, France.
\par
\noindent
$^2$ Groupe de Lyon: ENS Lyon, 46 all\'ee d'Italie, F-69364 Lyon Cedex
07,France.}\\
\end{center}

\vspace{20mm}

\centerline{ {\bf Abstract}}

\vspace{5mm}

\indent After a review of several methods designed to produce equivariant
cohomology classes, we apply one introduced by Berline, Getzler
and Vergne to get a family of representatives of the universal Thom
class of a vector bundle. Surprisingly, this family does not contain
the representative given by Matha\"{\i} and Quillen. However it
contains the particularly simple and symmetric representative of Harvey
and Lawson.

\vfill
\rightline{\LAP-A-574/96} \rightline{ 1996}

\newpage
\pagestyle{plain} \renewcommand{\thefootnote}{\arabic{footnote}} \newpage

\pagestyle{plain} \pagenumbering{arabic}

\section{Introduction.}

In a recent paper \cite{STW94} it has been shown how equivariant
cohomology is related to the so-called (cohomological) topological
models \cite{Bl92,BS88,BS91,OSB89,W88,WBS88}. In the same work, a way
to compute some representatives of equivariant cohomology classes
(i.e. observables of the corresponding topological model) was exhibited.

Here, we shall use this method in order to generate a  family of
representatives of the Thom class of a vector bundle depending on two
arbitrary functions. As we shall see, these
representatives are quite different from the Matha\"{\i}--Quillen
representative. They offer a good deal of
flexibility at the price of being slightly complicated. Special
choices allow to find a very special representative with remarkable
symmetry properties. However, its slow decrease at infinity makes it
necessary to consider a cohomology theory with coefficients with
sufficiently fast decrease (instead of compact). Some of these
representatives (in particular the most symmetric one) already appeared in
a quite different framework in the work of Harvey and Lawson on singular
connections \cite{hl93}, a fact we learned after this work was completed.  

This work is divided into three parts. In the first section we recall
basic facts about equivariant cohomology as well as the way
to compute representatives of equivariant cohomology classes. This
section parallels the explanations given in \cite{STW94}. The second
section is devoted to the Matha\"{\i}--Quillen representative of the
Thom class. Finally, the last section exhibits a large family of
representatives of the Thom class. 
\\

\section{Equivariant cohomology.} \label{sec:equicohom}

Let us consider the following setting~: ${\cal M}$ is a smooth manifold
and ${\cal G}$ a connected Lie group acting smoothly on ${\cal M}$. 
We would like to define a cohomology of the quotient space ${\cal M}/{\cal G}$ 
which coincides with the De Rham cohomology when this quotient is a smooth 
manifold but which also exists when it is not, i.e. when ${\cal G}$ acts with fixed
points. Equivariant cohomology solves this problem.

Let ${\cal M}$ be a smooth manifold and $\Om^* ({\cal M})$ the exterior
 algebra of differential forms on ${\cal M}$ endowed with the differential
$d_{{\cal M}}$. A Lie group ${\cal G}$ is assumed to be acting on ${\cal
M}$ as well as its Lie algebra, denoted $Lie{\cal G}$. For any $\lam 
\in Lie {\cal G}$ there is a vector field $\lam _{{\cal M}}$
representing the infinitesimal action of $\lam $ on ${\cal M}$.
This vector field $\lam _{{\cal M}}$ is usually called the
fundamental vector field associated with $\lam $. We shall denote
by $i_{{\cal M}}(\lam )=i_{{\cal M}}(\lam _{{\cal M}})$ and
$l_{{\cal M}}(\lam )=l_{{\cal M}}(\lam _{{\cal M}})=\left[
d_{{\cal M}},i_{{\cal M}}(\lam )\right] _{+}$ the contraction 
(or inner derivative) and Lie derivative acting on $\Om^* ({\cal M})$.
Let us recall that $i_{{\cal M}}(\lam )$ takes n-forms into (n-1)-forms
while $l_{{\cal M}}(\lam )$ acts on forms without changing the degree.
Elements of $\Om^* ({\cal M})$ which are annihilated by both $i_{{\cal
M}}(\lam )$ and $l_{{\cal M}}(\lam )$, for any $\lam \in Lie{\cal
G}$, are the so-called {\bf basic} elements of $\Om^* ({\cal M})$ for the
action of ${\cal G}$. As $d_{{\cal M}}$ maps basic elements into basic
elements, this leads to the definition of the {\bf basic
cohomology} of ${\cal M}$ for the action of ${\cal G}$ \cite{C50}.

We now consider the Weil algebra ${\cal W(\QTR{mathcal}{G})}$ of $Lie{\cal
G}$ \footnote{This is a harmless abuse of notation, but it is to be
remembered that equivariant cohomology deals only with the local
structure of ${\cal G}$.}. It is a graded differential algebra
generated by two $Lie{\cal G}$-valued indeterminates, the
``connection'' $\om $, of degree 1, and its 
``curvature'' $\Om $, of degree 2, such that~:  
\beq
\Om =d_{{\cal W}}\om +\dfrac 12\left[ \om ,\om \right]
\label{2.0}
\eeq
where $d_{{\cal W}}$ is the differential of ${\cal
W(\QTR{mathcal}{G})}$. Of course, one has the Bianchi identity~: 
\beq
d_{{\cal W}}\Om +\left[ \om ,\Om \right] =0 \label{2.1}
\eeq
There is an action $i_{{\cal W}}(\lam )$, $l_{{\cal
W}}(\lam )$ for $\lam \in Lie{\cal G}$~: 
\bea
i_{{\cal W}}(\lam )\om &=&\lam \qquad l_{{\cal W}}(\lam )\om 
=-[\lam ,\om ] \label{2.3} \\
i_{{\cal W}}(\lam )\Om &=&0 \qquad l_{{\cal W}}(\lam )\Om =-[\lam 
,\Om ] \label{2.4}
\eea
For instance, $\om $ may be a connection on a principal
${\cal G}$-bundle $\Pi $ and $\Om $ its curvature. In that case
$i_{{\cal W}}(\lam )$ and $l_{{\cal W}}(\lam )$ are generated by
the action of ${\cal G}$ on $\Pi $, and in this case
 ${\cal W(\QTR{mathcal}{G})}$ will be referred to as ${\cal W}_\Pi $. 

We now consider the graded differential algebra $(\Om^* ({\cal
M})\otimes {\cal W(\QTR{mathcal}{G})},d_{{\cal 
M}}+d_{{\cal W}})$, on which the operations $i_{{\cal M}}+i_{{\cal
W}}$ and $l_{{\cal M}}+l_{{\cal W}}(\lam )$ for any $\lam 
\in Lie{\cal G}$ are well-defined. There common
kernel is a graded differential subalgebra of $\Om^* ({\cal M})\otimes {\cal
W(\QTR{mathcal}{G})}$. By definition, the so-called {\bf 
equivariant cochains} are the elements of this subalgebra annihilated
by the differential $d_{{\cal M}}+d_{{\cal W}}$, leading to the
{\bf equivariant cohomology} of ${\cal M}$ for the action of ${\cal G}$~:
this is the so-called  {\bf Weil model} for equivariant cohomology.

Equivariant cohomology can be alternatively described in the so-called
{\bf intermediate model}, which was introduced in \cite{K93} and which will be
repeatedly used in the sequel. It is obtained from the Weil model via of
the following algebra isomorphism \footnote{See \cite{DV93} for a more general theorem.}~: 
\beq
x\longmapsto \exp \left\{ -i_{{\cal M}}(\lam )\right\} x \label{2.5}
\eeq
for any $x\in $ $\Om^* ({\cal M})\otimes {\cal W(\QTR{mathcal}{G})}$.
This isomorphism changes the original differential and operations on 
$\Om^* ({\cal M})\otimes {\cal W(\QTR{mathcal}{G})}$ by conjugation~: 
\bea
d_{{\cal M}}+d_{{\cal W}} &\longrightarrow &D_{int}^{}=d_{{\cal
M}}+d_{{\cal W}}+l_{{\cal M}}(\om )-i_{{\cal M}}(\Om )) \label{2.6} \\
(i_{{\cal M}}+i_{{\cal W}})(\lam ) &\longrightarrow &i_{{\cal W}}(\lam 
)=e^{-i_{{\cal M}}(\lam )}(i_{{\cal M}}+i_{{\cal W}})(\lam )e^{i_{
{\cal M}}(\lam )} \label{2.7} \\
(l_{{\cal M}}+l_{{\cal W}})(\lam ) &\longrightarrow &(l_{{\cal M}}+l_{
{\cal W}})(\lam )=e^{-i_{{\cal M}}(\lam )}(l_{{\cal M}}+l_{{\cal W}
})(\lam )e^{i_{{\cal M}}(\lam )} \label{2.8}
\eea

Finally, the so-called {\bf Cartan model} is obtained from the
intermediate model by putting $\om =0$ so that $D_{int}^2\mid
_{\om =0}$ vanishes when restricted to invariant cochains. This is
the most popular model, although many calculations are better
automatized in the intermediate model.

Another item which will be repeatedly used is ''Cartan's theorem 3'' \cite
{C50}~: let us assume that $(\Om^* ({\cal M}),d_{{\cal M}},i_{{\cal
M}},l_{{\cal M}})$ admits a ${\cal G}$-connection $\theta$ \footnote{
that is to say a $Lie{\cal G}$-valued 1-form on ${\cal M}$ such that $i_{
{\cal M}}(\lam )\theta =\lam $ and $l_{{\cal M}}(\lam )\theta
=-\left[ \lam ,\theta \right] $ for any $\lam \in Lie{\cal G}$.},
with curvature $\Theta $. Then any equivariant cohomology class of $\Om^* 
({\cal M})\otimes {\cal W(\QTR{mathcal}{G})}$ with representative $P(\om 
,\Om )$ gives rise canonically to a basic cohomology class of $\Om 
({\cal M})$ with representative $P(\theta ,\Theta )$. There is a
simple proof using the homotopy that expresses the triviality of the
cohomology of the Weil algebra \cite{MSZ85}. It follows from the
construction that the cohomology class of $P(\theta ,\Theta )$ does not
depend on $\theta $.

One convenient way to produce equivariant cohomology classes is as follows
\cite{BGV91}~: we consider an $H$-bundle ${\cal P}({\cal M},H)$ over {\cal M}
on which there exists an action of ${\cal G}$ which lifts the action of ${\cal G}$
on ${\cal M}$. In general, the Lie group $H$ has nothing to do with the Lie
group ${\cal G}$. As before, ${\cal P}({\cal M},H)$ is endowed with a
differential $d_{{\cal P}}$, a contraction $i_{{\cal P}}$ and a Lie
derivative $l_{{\cal P}}$.

Next, let $\Ga $ be a ${\cal G}$-invariant $H$-connection on ${\cal
P}({\cal M},H)$~: 
\beq
l_{{\cal P}}(\lam )\Ga =0\text{\quad ,\quad for any }\lam \in
Lie%
{\cal G} \label{2.9}
\eeq

The pull-back $\hat{\Ga }$ of $\Ga $ on $\Om^* ({\cal M})\otimes
{\cal W(\QTR{mathcal}{G})}$ is a 1-form on ${\cal P}({\cal M},H)$ and
a 0-form in ${\cal W(\QTR{mathcal}{G})}$. It follows that~: 
\beq
i_{{\cal W}}(\lam )\hat{\Ga }=0 \label{2.10}
\eeq
for any $\lam \in Lie{\cal G}$.

In $\Om^* ({\cal M})\otimes {\cal W(\QTR{mathcal}{G})}$, the equivariant
curvature of $\hat{\Ga }$ is defined by~: 
\beq
R_{int}^{eq}(\hat{\Ga },\om ,\Om )=D_{int}^{}\hat{\Ga }+\dfrac
12\left[ \hat{\Ga },\hat{\Ga }\right] \label{2.11}
\eeq
where $D_{int}^{}=d_{{\cal W}}+d_{{\cal P}}+l_{{\cal P}}(\om 
)-i_{{\cal P}}(\Om )$. Then, if $I_H$ is a symmetric invariant
polynomial on $LieH$, we consider the $H$-characteristic class
$I_{H,int}^{eq}(\hat{\Ga },\om ,\Om )=I_H\left(
R_{H,int}^{eq}(\hat{\Ga },\om ,\Om )\right) $. It is defined
on ${\cal M}$ and fulfills~: 
\bea
\left( d_{{\cal W}}+d_{{\cal M}}+l_{{\cal M}}(\om )-i_{{\cal M}}(\Om 
)\right) I_{H,int}^{eq}(\hat{\Ga },\om ,\Om ) &=&0 \label{2.12}
\\
i_{{\cal W}}(\lam )\ I_{H,int}^{eq}(\hat{\Ga },\om ,\Om ) &=&0
\label{2.13} \\
(l_{{\cal W}}+l_{{\cal M}})(\lam )\ I_{H,int}^{eq}(\hat{\Ga },\om 
,\Om ) &=&0 \label{2.14}
\eea
for any $\lam \in Lie{\cal G}$.

In the Weil model, the equivariant curvature is defined by~: 
\beq
R_W^{eq}(\hat{\Ga },\om ,\Om )=(d_{{\cal W}}+d_{{\cal
P}})\hat{\Ga }+\dfrac 12\left[ \hat{\Ga }+i_{{\cal P}}(\Om 
)\hat{\Ga },\hat{\Ga } +i_{{\cal P}}(\Om )\hat{\Ga }\right]
\label{2.15} 
\eeq

We may similarly consider~: 
\beq
I_{H,W}^{eq}(\hat{\Ga },\om ,\Om )=I_H\left(
R_{H,W}^{eq}(\hat{\Ga },\om ,\Om )\right) =e^{-i_{{\cal M}}(\lam 
)}I_{H,int}^{eq}(\hat{\Ga },\om ,\Om ) \label{2.16}
\eeq
which fulfills~: 
\bea
(d_{{\cal W}}+d_{{\cal M}})\ I_{H,W}^{eq}(\hat{\Ga },\om ,\Om )
&=&0 \label{2.17} \\
(i_{{\cal W}}+i_{{\cal M}})(\lam )\ I_{H,W}^{eq}(\hat{\Ga },\om 
,\Om ) &=&0 \label{2.18} \\
(l_{{\cal W}}+l_{{\cal M}})(\lam )\ I_{H,W}^{eq}(\hat{\Ga },\om 
,\Om ) &=&0 \label{2.19}
\eea
for any $\lam \in Lie{\cal G}$.

Finally, if ${\cal M}$ admits a ${\cal G}$-connection $\theta $ with
curvature $\Theta $, we can apply ''Cartan's theorem 3'', and substitute
$\theta $ and $\Theta $ instead of $\om $ and $\Om $ in
$I_{H,W}^{eq}(\hat{\Ga },\om ,\Om )$, so that~: 
\bea
d_{{\cal M}}\ I_{H,W}^{eq}(\hat{\Ga },\theta ,\Theta ) &=&0 
\label{2.20} \\
i_{{\cal M}}(\lam )\ I_{H,W}^{eq}(\hat{\Ga },\theta ,\Theta ) &=&0
\label{2.21} \\
l_{{\cal M}}(\lam )\ I_{H,W}^{eq}(\hat{\Ga },\theta ,\Theta ) &=&0
\label{2.22}
\eea
for any $\lam \in Lie{\cal G}$.

By standard arguments, these cohomology classes do not depend either on $%
\hat{\Ga }$ or on $\theta $.

\section{The Thom Class of a Vector Bundles~: the Matha\"{\i}--Quillen
strategy \protect\cite{MQ86}.}

Let $V$ be a real oriented Euclidean vector space of dimension $n=2d$ with
scalar product $(\ ,\ )_V$. On $V$, we choose a canonical basis
$\left\{ \vec{e}_k\right\} $ orthonormal with respect to $(\ ,\ )_V$~: 
\beq
(\vec{e}_i,\vec{e}_j)_V=\delta _{ij} \label{3.1}
\eeq
Any vector on $V$ can be  decomposed as~: 
\beq
\vec{v}=v^k\vec{e}_k \label{3.2}
\eeq
Such a decomposition gives a coordinates system $(v^k)$ on $V$,
turning $V$ into a manifold. Due to the linear space structure of $V$, only
$GL(n,{\Bbb R})$ transformations define allowed coordinate changes. The
group of isometries of $V$, with respect to $(\ ,\ )_V$, is
$SO(n)\subset GL(n,{\Bbb R})$, with Lie algebra $so(n)$ and Weil
algebra ${\cal W}(SO(n))$. Finally, we endow $V$ and ${\cal W}(SO(n))$
with the standard differential operations $d_V$, $i_V$, $l_V$, $d_{{\cal W}}$,
$i_{{\cal W}}$ and $l_{{\cal W}}$.

Now, let $E({\cal M},V)$ be a vector bundle over a smooth manifold ${\cal
M}$ with typical fiber $V$, equipped with differential operations~: $d_E$,
$i_E$ and $l_E$,. We denote $\Om _{rdv}^n(E)$ the space of n-forms on $E$
whose restriction to each fiber of $E$ is rapidly decreasing. The
corresponding cohomology space is written $H_{rdv}^n(E)$. The Thom
Class of $E$ is the element T$(E)$ of $H_{rdv}^n(E)$ such that~: 
\beq
\int_VT(E)=1 \label{3.3}
\eeq
which means that integration of T$(E)$ along the fiber produces
the constant function 1 on ${\cal M}$.

Actually, following Matha\"{\i} and Quillen \cite{MQ86}, we would like to
exhibit a representative of T$(E)$ in the form of an integral representation.
Then, we consider $V^{*}$, the
dual space of $V$, equipped with the scalar product $(\ ,\ )_{V^{*}}$,
dual to $(\ ,\ )_V$ on $V$. Moreover, we introduce coordinates $(\varpi _k)$
for the Grassmann algebra $\Lam V^{*}$ of $V^{*}$ together with the
differential operations $\delta $, $I$ and $L$, dual to those on
$V$.

We take as structure equations~:
\beq 
\ba{rcl} \label{eq:struceq}
s^{top}v^k &=&\Psi ^k+L^{top}(\om )v^k  \\
s^{top}\Psi ^k &=&-L^{top}(\Om )v^k+L^{top}(\om )\Psi ^k 
 \\
s^{top}\varpi _k &=&b_k+L^{top}(\om )\varpi _k  \\
s^{top}b_k &=&-L^{top}(\Om )\varpi _k+L^{top}(\om )b_k 
\\
s^{top}\om &=&\Om -\dfrac 12\left[ \om ,\om \right] 
 \\
s^{top}\Om &=&-\left[ \om ,\Om \right] 
\ea
\eeq
with~: 
\bea
s^{top} &=&d_{{\cal W}}+(d_V+\delta )+(l_V+L)(\om )-(i_V+I)(\Om )
\label{3.10} \\
\Psi ^k &=&d_Vv^k\equiv \Psi _{int}^k \label{3.11}
\eea
in the intermediate model, and~: 
\bea
s^{top} &=&d_{{\cal W}}+d_V+\delta \label{3.12} \\
\Psi ^k &=&(d_V-L^{top}(\om ))v^k\equiv \Psi _W^k \label{3.13}
\eea
in the Weil model, while~: 
\beq
L^{top}=l_V+L \label{3.14}
\eeq
in any model.

The null section $s_0$ of $E({\cal M},V)$ that sends any point of ${\cal
M}$ into the null vector, diffeomorphically maps ${\cal M}$ into $s_0({\cal M}
)\subset E$. Then, the Thom Class T$(E)$ of $E$ is nothing but the
Poincar\'{e} dual of $s_0({\cal M})$ in $E$ \cite{BT82}, and the Dirac
form on $E$~: 
\beq
\delta (\vec{v})dv^1\wedge ...\wedge dv^n \label{3.15}
\eeq
represents the Poincar\'{e} dual of $s_0({\cal M})$ in $E$. This
form can be written as a Fourier transform~: 
\beq
\frac 1{(2\pi )^n}\int dbd\varpi \exp i\left\{ b.\vec{v}+\varpi .\Psi
\right\} =\frac 1{(2\pi )^n}\int dbd\varpi \exp i\left\{ b_kv^k+\varpi
_k\Psi ^k\right\} \label{3.16}
\eeq
From the structure equations (\ref{eq:struceq}), we deduce~: 
\beq
b.\vec{v}+\varpi .\Psi =s^{top}(\varpi .\vec{v}) \label{3.17}
\eeq
However, we can consider a smoother representative, with a
gaussian behavior for instance. That means that we must insert a term of
the form~: 
\beq
i(b,b)_{V^{*}} \label{3.18}
\eeq
into (\ref{3.17}). Now, we can try to write the new argument as a
$s^{top}$-exact term~: 
\beq
s^{top}(\varpi .\vec{v}+i(\varpi ,b)_{V_{*}})=b.\vec{v}+\varpi .\Psi
+i(b,b)_{V^{*}}-i(L^{top}(\Om )\varpi ,\varpi )_{V^{*}} \label{3.19}
\eeq
so that we are led to define~: 
\beq
U=\frac 1{(2\pi )^n}\int dbd\varpi \exp i\left\{ s^{top}(\varpi .\vec{v}%
+i(\varpi ,b)_{V_{*}})\right\} \label{3.20}
\eeq
Note that U is an element of ${\cal W}(SO(n))\otimes \Om^* (V)$.

In order to prove that U maps into a representative of T$(E)$, let us proceed in
the intermediate model where we write $U_{int}$ instead of $U$. Then,
since in (\ref{3.20}) $\om \in {\cal W}(SO(n))$ does not appear, we
immediately conclude that $U_{int}$ does not explicitly depend on
$\om $, that is to say~: 
\beq
\forall \lam \in so(n)\qquad i_{{\cal W}}(\lam )U_{int}=0 
\label{3.21}
\eeq
which express the basicity condition within the intermediate
model. Now, there remains to show that $U_{int}$ is closed with respect to
$D_{int}=d_{{\cal W}}+d_V+l_V(\om )-i_V(\Om )$. Indeed~: 
\bea
D_{int}U_{int} &=&\frac 1{(2\pi )^n}D_{int}\int dbd\varpi \exp \left\{
i.s^{top}(\varpi .\vec{v}+i(\varpi ,b)_{V_{*}})\right\} \label{3.22} \\
&=&\frac 1{(2\pi )^n}\int dbd\varpi \left( s^{top}-D_{V^{*}}\right) \exp
\left\{ i.s^{top}(\varpi .\vec{v}+i(\varpi ,b)_{V_{*}})\right\} 
\label{3.23}
\eea
where $D_{V^{*}}=\delta +L(\om )-I(\Om )$. Hence~: 
\beq
D_{int}U_{int}=-\frac 1{(2\pi )^n}\int dbd\varpi \left[ D_{V^{*}}\exp
\left\{ i.s^{top}(\varpi .\vec{v}+i(\varpi ,b)_{V_{*}})\right\} \right]
\label{3.24}
\eeq
Now, from the structure equations (\ref{eq:struceq}), we get~: 
\bea
D_{V^{*}} &=&\left( b_k+L^{top}(\om )\varpi _k\right) \dfrac \partial
{\partial \varpi _k}+\left( -L^{top}(\Om )\varpi _k+L^{top}(\om 
)b_k\right) \dfrac \partial {\partial b_k} \label{3.25} \\
&=&\left( L^{top}(\om )\varpi _k\dfrac \partial {\partial \varpi
_k}+L^{top}(\om )b_k\dfrac \partial {\partial b_k}\right) +\left(
b_k\dfrac \partial {\partial \varpi _k}-L^{top}(\Om )\varpi _k\dfrac
\partial {\partial b_k}\right) \label{3.26}
\eea
The first term in $D_{V^{*}}$ corresponds to an
$so(n)$-transformation. Due to the $so(n)$-invariance of the
measure $dbd\varpi $, it does not contribute to (\ref{3.24}). The last
term in (\ref{3.26}) vanishes upon integration by parts. Then~: 
\beq
D_{int}U_{int}=0 \label{3.27}
\eeq
Finally, combining equations (\ref{3.21}) and (\ref{3.27}), we deduce that~:

\beq
\forall \lam \in so(n)\qquad (l_{{\cal W}}+l_V)(\lam )U_{int}=0
\label{3.28}
\eeq
and conclude that $U_{int}$ is a representative in ${\cal W}
(SO(n))\otimes \Om^* (V)$ of the Thom Class of $E({\cal M},V)$. The
corresponding representative in the Weil model is obtained by setting~: 
\bea
s^{top} &=&d_{{\cal W}}+d_V+\delta \label{3.29} \\
\Psi ^k &=&(d_V-L^{top}(\om ))v^k\equiv \Psi _W^k \label{3.30}
\eea
within equation (\ref{3.20}).

\indent

Actually, it can be easily shown that Fourier transform (denoted ${\cal F}$) commutes with
equivariant differential operations. More precisely :

\bea
{\cal F}\left[ \left( d_{{\cal W}}+\delta +L(\om )-I(\Om )\right) \Phi
\right] &=&\left( d_{{\cal W}}+d_V+l_V(\om )-i_V(\Om )\right)
{\cal F}\left[
\Phi \right] \label{3.31} \\
{\cal F}\left[ i_{{\cal W}}(\lam )\Phi \right] &=&i_{{\cal W}}(\lam 
){\cal F}\left[ \Phi \right] \\
{\cal F}\left[ (l_{{\cal W}}+L)(\lam )\Phi \right] &=&(l_{{\cal W}%
}+l_V)(\lam ){\cal F}\left[ \Phi \right]
\eea
in the intermediate model. The same holds in the Weil model with
suitable differentials. Let us point out that this mainly relies on the 
identity $b.\vec{v}+\varpi .\Psi =s^{top}(\varpi .\vec{v})$.

Then, since $\phi =(b,b)_{V^{*}}+(L^{top}(\Om )\varpi ,\varpi
)_{V^{*}}$ is equivariant, it is straightforward to find that its
Fourier transform is also equivariant. This simple remark allows to
construct representatives of equivariant cohomology classes using
Fourier transform of functions of $\phi$.

\indent

Finally, we can consider a principal $SO(n)$-bundle $P$ over ${\cal M}$.
It is well known that $P\times _{SO(n)}V$ is a vector bundle isomorphic to
$E$, and $P\times V$ is called the principal $SO(n)$-bundle associated with
$E({\cal M},V)$. Hence, as an n-form on $E$, any representative of the Thom
Class T$(E)$ of $E$ comes from a closed $SO(n)$-basic $n$-form on the
associated bundle $P\times V$ of $E$. In order to produce such a representative of
T$(E)$, we use Cartan's Theorem 3, that is to say we replace $(\om 
,\Om ) $ (in the representative $U$) by $(\theta ,\Theta )$, a
connection and its curvature on $P({\cal M},SO(n))$.

\section{Construction of Representatives of Thom Class of Vector Bundles~:
the Berline--Getzler--Vergne strategy \protect\cite{BGV91}.}

In this section, we shall use the strategy explained in section 2) in
order to produce representatives of T$(E)$.

To begin with, we are going to turn $V$ into a Riemannian manifold, i.e. a
manifold $V$ with a metric. The tangent bundle of $V$, denoted by $TV$, is
obviously isomorphic to $V\times V$. The only $SO(n)$-invariants formed
with $\vec{v}$ and $d\vec{v}$ are the three scalar products,
so that the general $SO(n)$-invariant metric on $V$ is~:

\beq
ds^2(\vec{v})=e^{\vf }\left( (dv^i)^2+\si (v^idv^i)^2\right)
\label{4.1}
\eeq
where $\vf $ and $\si $ are smooth functions of
$t=(\vec{v},\vec{v})_V$ only. The above expression is positive
definite if and only if $1+\si (t)t > 0$ for $t \geq 0$. One can
assume if convenient that the metric is asymptotically flat (i.e. that
the curvature vanishes at infinity).

We can consider the principal $GL(n,{\Bbb R})$-bundle associated with
$TV$, i.e. the frame bundle $R(V)$ of $V$. It is made of the points
$(\vec{v},b_{\vec{v}})$ where $b_{\vec{v}}$ is a frame (i.e. a basis)
at $\vec{v}$. Coordinates for $b_{\vec{v}}$ are defined as 
follows. We denoted by $(\partial _k)$ the natural basis of $T_{
\vec{v}}V$ defined by the canonical coordinates $(v^k)$ of $V$~: $\partial
_k=\frac \partial {\partial v^k}$. Then, the coordinates of $b_{\vec{v}}$
are the components $b_k^j$ of the decomposition of $b_{\vec{v}}$ with
respect to the natural basis $(\partial _k)$~: 
\beq
b_{\vec{v}\ k}=b_k^j\partial _j \label{4.2}
\eeq
with $b_{\vec{v}\ k}$ the $k$-th frame vector of the frame
$b_{\vec{v}}$. The isometry group of $(V,(\ ,\ )_V)$, namely $SO(n)$,
acts both on elements of $V$ and on frames, that is to say on $R(V)$.
This goes as follows. For any $\Phi \in SO(n)$~: 
\beq
\Phi ^k(\vec{v})=\Phi _m^kv^m \label{4.3}
\eeq

At the infinitesimal level, if we write $\Phi _m^k=\delta _m^k+\vf
_m^k$, we get~: 
\beq
\Phi ^k(\vec{v})=\left( \delta _m^k+\vf _m^k\right) v^m=v_{}^k+\vf
_m^kv^m=v_{}^k+\xi ^k \label{4.4}
\eeq
where $\xi ^k=$ $\vf _m^kv^m$ defines a vector field on $V$,
the so-called fundamental vector field associated with the action of $
\vf \in so(n)$~: 
\beq
\xi =\xi ^k\dfrac \partial {\partial v^k} \label{4.5}
\eeq
The natural action of $\Phi \in SO(n)$ on $T_{\vec{v}}V$, 
is given by the so-called differential of $\Phi $ at $\vec{v}$,
$d_{\vec{v}}\Phi :T_{\vec{v}}V\rightarrow T_{\Phi (\vec{v})}V$~: 
\beq
\forall X_{\vec{v}}\in T_{\vec{v}}V\quad ,\quad \forall f\in C^\infty
(V)\quad \quad d_{\vec{v}}\Phi X_{\vec{v}}\ (f)=X_{\vec{v}}(f\circ \Phi )
\label{4.6}
\eeq
Applying this definition to the frame vectors $b_{\vec{v}\ k}$,
one gets~: 
\beq
\tilde b_j^i=b_j^m\left( \partial _m\Phi ^i(\vec{v})\right) =b_j^m\Phi
_m^i
\label{4.7}
\eeq
where $\tilde b_j^i$ are the coordinates of the transformed frame
at $\Phi (\vec{v})$. At the infinitesimal level, for $\vf \in so(n)$~:

\beq
\tilde b_j^i=b_j^m\left( \delta _m^i+\vf _m^i\right)
=b_j^i+b_j^m\vf _m^i=b_j^i+b_j^m\vf _m^i=b_j^i+\Xi _j^i 
\label{4.8}
\eeq
Combining equations (\ref{4.4}) and (\ref{4.8}), we deduce that the fundamental
vector field associated with the action of $\vf \in so(n)$ on $R(V)$
reads~: 
\beq
\lam _R=\xi ^k\dfrac \partial {\partial v^k}+\Xi _q^p\dfrac \partial
{\partial b_q^p}=\vf _m^kv^m\dfrac \partial {\partial
v^k}+b_q^m\vf
_m^p\dfrac \partial {\partial b_q^p} \label{4.9}
\eeq

Now, let $P({\cal M},SO(n))$ be some principal $SO(n)$-bundle over a
smooth manifold ${\cal M}$. It is well known that there is a vector
bundle over ${\cal M}$ associated with $P$ for the action of $SO(n)$
on $V$. The group $SO(n)$ acts on the right on $P$ and on the left on
$V$. We first define a right-action of $SO(n) $ on $P\times V$ by setting~: 
\beq
(p,\vec{v})^\Phi =(p.\Phi ,\Phi ^{-1}(\vec{v})) \label{4.10}
\eeq
so that, the fundamental vector field representing the action of
$\vf \in so(n)$ on $P\times V$ reads~: 
\beq
\lam _{P\times V}=\lam _P-\xi ^k\dfrac \partial {\partial
v^k}=\lam 
_P-\vf _m^kv^m\dfrac \partial {\partial v^k} \label{4.11}
\eeq
where $\lam _P$ is the fundamental vector field representing
the action of $\vf $ on $P$.

Finally, the action of any $\vf \in so(n)$ on the $GL(n,{\Bbb
R})$-principal bundle $P\times R(V)$, is given by following
fundamental vector field~: 
\beq
\lam =\lam _P-\lam _R \label{4.12}
\eeq
with $\lam _R$ defined in equation (\ref{4.9}).

In the following, $V$, $R(V)$ and $P$ are equipped with the following
differential operations~: $d_V$, $d_R$, $d_P$, $i_V$, $i_R$, $i_P$, $l_V$,
$l_R$ and $l_P$, respectively exterior differentials, inner products and
Lie derivatives.

Now, since we are looking for representatives of equivariant cohomology
classes, we can mimic the construction made in \cite{STW94} in the case of
two-dimensional Gravity. We first look for a $GL(n,{\Bbb R})$-connection
on ${\cal P}(P\times V,GL(n,{\Bbb R}))=P\times R(V)$ invariant under the
action of $SO(n)$. If we notice that, by construction, the metric {\bf
g} on $V$ is $SO(n)$-invariant, we can consider the Levi-Cevita
connection $^{LC}\Ga $ associated with {\bf g}. Due to the
$SO(n)$-invariance of {\bf g}, $^{LC}\Ga $ is an $SO(n)$-invariant
connection. More precisely, the lift of $^{LC}\Ga $ into a
connection one form $\Ga $ on $R(V)$ according to~: 
\beq
\Ga =b^{-1}(^{LC}\Ga )b+b^{-1}d_Rb \label{4.13}
\eeq
is invariant under the action of $SO(n)$. The fundamental vector
field for the action of $so(n)$ was given before, so that~: 
\beq
\left( i_{{\cal P}}(\lam )\Ga \right) _\tau ^\si =(b^{-1})_\nu
^\si \left( -\ ^{LC}D_\mu \xi ^\nu \right) b_\tau ^\mu \label{4.14}
\eeq
where $i_{{\cal P}}(\lam )=(i_P+i_R)(\lam )$, and~: 
\beq
l_{{\cal P}}(\lam )\Ga =0 \label{4.15}
\eeq
with $l_{{\cal P}}(\lam )=(l_P+l_R)(\lam )$.

The next step is to consider the Weil algebra ${\cal W}(SO(n))$ of
$so(n)$. The relevant formul\ae\ were given in section
\ref{sec:equicohom}. We recall that  the equivariant curvature of $\Ga
$ in the intermediate model is~:

\beq
R_{int}^{eq}(\Ga ,\om ,\Om )=\left( d_{{\cal W}}+d_{{\cal
P}}+l_{{\cal P}}(\om )-i_{{\cal P}}(\Om )\right) \Ga +\dfrac
12\left[\Ga ,\Ga \right] \label{4.20}
\eeq
while the corresponding curvature in the Weil model is obtained as~: 
\beq
R_W^{eq}(\Ga ,\om ,\Om )=e^{i_{{\cal P}}(\om 
)}R_{int}^{eq}(\Ga ,\om ,\Om ) \label{4.21} 
\eeq
which gives~: 
\beq
R_W^{eq}(\Ga ,\om ,\Om )=R(\Ga )+i_{{\cal P}}(\om )R(\Ga 
)+\dfrac 12i_{{\cal P}}(\om )i_{{\cal P}}(\om )R(\Ga )-i_{{\cal
P}}(\Om )\Ga \label{4.22}
\eeq
The Weil equivariant Euler class is defined by~:

\beq
E_W^{eq}=\dfrac{\ve ^{\mu _1\rho _1...\mu _d\rho _d}}{\sqrt{{\bf
g}}}g_{\rho _1\upsilon _1}...g_{\rho _d\upsilon _d}\left( R_W^{eq}\right)
_{\mu_1}^{\upsilon _1}\wedge ...\wedge \left( R_W^{eq}\right) _{\mu
d}^{\upsilon_d} \label{4.23}
\eeq
which after normalization gives rise to a representative of T$(E)$ in
$P({\cal M},SO(n))\times V$. 

It is now time to use the explicit form of the metric to get a formula
for the Thom class. Surprisingly, we shall see there is no choice of
metric that allows to recover the Matha\"{\i}--Quillen representative of
T$(E)$. From now on, the computations, if painful, are
straightforward. We use the intermediate model so that $d_V v^i\equiv
\Psi^i$. As (\ref{4.13}) looks formally like a change of
coordinates in the fiber, we know that its effect on curvature will
be a simple conjugation which disappears completely on the Thom class.
So we can forget it in the computation. From (\ref{4.1}) we find that
the metric is~: 

\beq
g_{ij}=e^{\vf }(\delta_{ij}+\si v_iv_j)
\eeq

Our notations need some comment~: we start with global coordinates
$v^i$ on $V$, so the exponent $i$ is not a tensor component but just a
label. The metric is expressed with respect to this particular
coordinate system. However it is convenient to deal consistently with
formal lower and upper indices in the Einstein summation
convention. So we define $v_i \equiv v^i$ and
$\delta_{ij}\equiv\delta^{ij}\equiv\delta^i_j\equiv\delta_i^j=1$ if
$i=j$ and $0$ else. For instance we use the notation $v_i$ and
$\delta_{ij}$ in $g_{ij}$ and we write $t=v_iv^i$. This becomes
slightly less formal if we restrict the diffeomorphism group of $V$ to
linear orthogonal transformations. 

A simple computation shows that the inverse metric is

\beq
g^{ij}=e^{-\vf }(\delta^{ij}+\tsi v^iv^j)
\eeq
where $\tsi$ is defined by $(1+t\tsi )(1+t\si )=1$.

First, we need a formula for the connection and curvature. The fact
that $\vf$ and $\si$ depend only on $t$ leads to many simplifications
in the computation. We use dots for derivatives with respect to $t$. 

With the expression of $\bf g$, we get for the
connection \footnote{Remember that $\Ga ^k_{ij}\equiv
\frac{1}{2}g^{kl}(\partial_i g_{lj} + \partial_j 
g_{il} - \partial_l g_{ij})$.}~:

\beq
\Ga ^k_{ij}=(1+t\tsi)v^k[(\si-\dot{\vf})\delta_{ij} +
(\dot{\si}-\si\dot{\vf}) v_iv_j]+\dot{\vf}(v_i\delta_j^k+v_j\delta_i^k).
\eeq
so that the connection matrix is
\beq 
\Ga _i^j\equiv \Psi^k\Ga _{ik}^j=Av^jv_iv_k\Psi^k
+Bv^j\Psi_i+C(v_i\Psi^j+\delta^j_iv_k\Psi^k)
\eeq
where we have set
\beq
\left\{ \ba{l} A=(1+t\tsi)(\dot{\si}-C\si) \\ B=(1+t\tsi)(\si -C) \\
C=\dot{\vf}\ea \right. . 
\eeq

The curvature matrix is given by 

\beq
R_i^j \equiv d\Ga _i^j+\Ga ^j_k \wedge \Ga ^k_i.
\eeq

A tedious computation leads to

\beq
R^{ij}\equiv g^{ik}R^j_k=e^{-\vf}(1+t\tsi) \left[ M_{11} \Psi^i\Psi^j
+ M_{12} \left(v^i(v_k\Psi^k)\Psi^j-v^j(v_k\Psi^k)\Psi^i\right)\right]
\eeq
where
\beq
\left\{ \ba{l} M_{11}=C^2t+2C-\si \\
M_{12}=2\dot{C}-C^2-\dot{\si }+(1+t\tsi)(\si -C)\dot{t\si} \ea \right. . 
\eeq

To get the full equivariant curvature, we need the part involving
$\Om $. In accordance with our convention on indices, we define
$\Om ^{ij}\equiv \Om _i^j$. By definition $\Om ^{ij}$ is antisymmetric.
According to formul\ae\ (\ref{4.14},\ref{4.20}), the part of the
equivariant curvature containing $\Om $
is the covariant derivative of $\Om _k^jv^k$, the $so(n)$ vector
field associated to $\Om $. Consequently

\beq
\left(-i_{{\cal P}}(\Om ) \Ga \right)_k^j=\Om _k^j+\Om _m^lv^m\Ga ^j_{lk}.
\eeq

The antisymmetry of $\Om $ leads to further simplifications. 
The outcome is~:

\beq
g^{ik}\left(\Om _k^j+\Om _m^lv^m\Ga ^j_{lk}\right)=e^{-\vf} (1+t\tsi)\left[
M_{21} \Om ^{ij} + M_{22} \left(v^iv_k\Om ^{kj}-v^jv_k\Om ^{ki}
\right)\right]
\eeq 
where
\beq
\left\{ \ba{l} M_{21}=1+t\si \\
M_{22}=C-\si \ea \right. . 
\eeq

We note the striking similarity between the two contributions.
If we define a two-by-two matrix $N^{ij}$ by

\beq
N^{ij}=\left(\ba{cc} \Psi^i\Psi^j & \Om ^{ij} \\
v^i(v_k\Psi^k)\Psi^j-v^j(v_k\Psi^k)\Psi^i & v^iv_k\Om ^{kj}-v^jv_k\Om
^{ki}\ea \right) 
\eeq
the equivariant curvature can be written as a trace ~:
\beq
(R_{int}^{eq})^{ij}=e^{-\vf}(1+t\tsi)\Tr  MN^{ij}
\eeq

The equivariant Euler class is
\beq
E_{int}^{eq}=2^{n/2}\sqrt{\bf g}\, \Pf (R_{int}^{eq})^{ij}.
\eeq
with the usual definition of the Pfaffian. Note that $\sqrt{\bf
g}=e^{n\vf/2}(1+t\si)^{1/2}$. 

This is the explicit formula for the universal Thom class that we
were after. It involves two arbitrary functions of $t$, $\vf$ and
$\si$ (with the mild restriction $1+t\si > 0$) which may be localized
at will thus so leaving a fair amount of flexibility. 

The first comment to make is that apparently the above representative,
which is of course $so(n)$ invariant when $so(n)$ acts on $V$, $\Om $ and
$\Psi$ at the same time, is not invariant when $so(n)$ acts only on
$V$. To state it more simply, the $V$ dependence of the Thom class is
not only through $t$. This is to be contrasted with the
Matha\"{\i}--Quillen representative.

Let us deal with a special case first. When $n =2$, it is easy
to see that 

\beq
\ve _{ij} N^{ij} =2\left(\ba{cc} \Psi^1\Psi^2 & \Om ^{12} \\
t\Psi^1\Psi^2 & t\Om ^{12}\ea \right) 
\eeq
so we have some hope to recover the Matha\"{\i}--Quillen formula
as a special case. After some manipulations one finds

\beq \label{eq:Ed=2}
E_{int}^{eq}=4\dot{F}\Psi^1\Psi^2+2F\Om ^{12}
\eeq
where

\beq \label{eq:F}
F \equiv \frac{1+tC}{(1+t\si)^{1/2}}.
\eeq
So the Thom class depends only on one arbitrary function of $t$,
namely $F$, which can easily be adjusted to recover the
Matha\"{\i}--Quillen representative. The correct choice is 
$F=-{\textstyle{1 \over {(4\pi )}}}\exp (-{t \mathord{\left/
{\vphantom {t {4)}}} \right. \kern-\nulldelimiterspace} {4)}}$.

When $n > 2$ the situation is more complicated. We shall use a trick
to see how much the symmetry of the $so(n)$ action on $V$ is broken.

The first observation is that under a similarity, the Pfaffian has a
simple behavior~: if $A$ is a square antisymmetric matrix and $S$
an arbitrary square matrix of the same size with transpose $S^t$,
$S^tAS$ is again antisymmetric and $\Pf S^tAS =\De S \; \Pf A$. The square
of this equation just follows from the 
multiplicative property of the determinant, and the sign is fixed
by the case when $S$ is the identity matrix. So if we can find a
matrix $S^{ij}$ (independent of $\Psi$ and $\Om$) such
that $S^tR_{int}^{eq}S$ simplifies, we shall end with a simpler formula
for the Thom class. 

Define a symmetric matrix $S(D)$ of parameter $D$ by

\beq
S(D)^i_j=\delta^i_j+Dv^iv_j.
\eeq
This matrix is easily diagonalized : the vectors orthogonal to $v^i$
are left invariant and $v^i$ is multiplied by $1+tD$. So

\beq \De S(D)=1+tD \eeq
and 

\beq S(D)S(E)=S(D+E+tDE). \eeq

Moreover, if $A^{ij}$ is any antisymmetric matrix, 

\beq 
(S(D)AS(D))^{ij}= A^{ij}+D(v^iv_kA^{kj}-v^jv_kA^{ki}). 
\eeq

We apply this identity to the four antisymmetric objects building
the two-by-two matrix $N^{ij}$ to get

\beq S(D)NS(D)=\bar{S}(D)N \label{ij22} \eeq
where $\bar{S}(D)$ is the  two-by-two matrix 
\beq  \left(\ba{cc} 1 & D \\ 0 & 1+tD \ea \right) .\eeq 

In equation (\ref{ij22}), the left-hand side
involves a product of $n \times n$ matrices, and the $2 \times 2$
indices are spectators whereas on the right-hand side the opposite
occurs.    

So we can write

\beq
S(D)R_{int}^{eq}S(D)=e^{-\vf}(1+t\tsi)\Tr M(D)N
\eeq
with $M(D)\equiv M\bar{S}(D)$, and
the Thom class is

\beq
E_{int}^{eq}=(1+t\si)^{(1-n)/2}(1+tD)^{-1}\Pf (\Tr M(D)N).
\eeq

We can choose $D$ to simplify the expression of $E_{int}^{eq}$.

First we take $D=D_1$ where $M(D_1)=\left(\ba{cc} * & 0 \\ {*} & * \ea
\right)$. 
This makes it easy to compute the term in $E_{int}^{eq}$ that does
not involve $\Om$. The outcome is

\beq
E_{int}^{eq}=n!(1+t\si)^{(1-n)/2}(1+t\frac{M_{12}}{M_{11}})
M_{11}^{n/2} \Psi^1\cdots\Psi^n + \mbox{\rm terms involving } \Om .
\eeq
One can check that this is compatible with (\ref{eq:Ed=2}) for $n=2$.

Second, we take $D=D_2$ where $M(D_2)=\left(\ba{cc} * & * \\ {*} & 0 \ea
\right)$. This makes it easy to compute the term in $E_{int}^{eq}$
that does not involve $\Psi $. The outcome is

\beq
E_{int}^{eq}=2^{n/2}(1+t\si)^{(1-n)/2}(1+t\frac{M_{22}}{M_{21}})
M_{21}^{n/2} \Pf \Om  + \mbox{\rm terms involving } \Psi
\eeq
a result which is again compatible with (\ref{eq:Ed=2}) for $n=2$.

Those two terms in $E_{int}^{eq}$ automatically depend only on $t$.
On the other hand, the other terms are not scalars for the $so(n)$
action on $V$. To see this we keep $D=D_2$, set
$A^{ij}=M_{11}\Psi^i\Psi^j+M_{21} \Om ^{ij}$ and
$B^{ij}=v^i(v_k\Psi^k)\Psi^j-v^j(v_k\Psi^k)\Psi^i$. Using the fact that
$x_k\Psi ^{k}$ squares to $0$ we get

\beq \label{eq:eulerclass}
\ba{rcl}
E_{int}^{eq}& =& \\ & & \hspace{-.7cm}(1+t\si)^{(1-n)/2}\left(
2^{n/2}(1+t\frac{M_{22}}{M_{21}})\Pf A - n/2 \frac{\De
M}{M_{21}}\ve_{i_1j_1\cdots i_nj_n}B^{i_1j_1}A^{i_2j_2}\cdots
A^{i_nj_n}\right).
\ea  
\eeq

As $\Om$ and $\Psi$ are independent families of indeterminates, the
matrix elements of $A^{ij}$ are independent of each other (except for
antisymmetry) and of the matrix elements of $B^{ij}$. So in the expansion
of 

\beq
\ve_{i_1j_1\cdots i_nj_n}B^{i_1j_1}A^{i_2j_2}\cdots A^{i_nj_n}
\eeq
no compensation can occur between $A$-factors and
$B$-factors or between different $B$-factors. Moreover $B$-factors
contain the full non $so(n)$ invariant part of the $V$ dependence
of the Thom class. So we have the following three possibilities. Either $\De
M$ is $0$, or $B^{ij}$ is invariant for the action of $so(n)$ on $V$,
or the representative of the Thom class is not invariant for the
action of $so(n)$ on $V$. The first term of the alternative depends on
our choice of $\vf$ and $\si$. The second is easily checked to occur
if and only if $n=2$, a case we have already treated.

So finally, we have shown that if $n > 2$ the representative of the
Thom class is invariant for the $so(n)$ action on $V$ if and only if
$\De M =0$. 

We shall now see that despite the fact that apparently our
representative of the Thom class depends on two arbitrary functions,
the single condition $\De M =0$ fixes it completely. This can be seen as
a manifestation of the topological character of the Thom class. We
shall also see that the representative we end up with is not the
Matha\"{\i}--Quillen representative.

From now on, we set $\De M =0$. Explicit computation shows that this
equation has a first integral. Namely $\De M =0$ is equivalent to

\beq \label{eq:detm=0}
(1+tC)^3\frac{d}{dt}
\left(\frac{\si}{(1+tC)^2}+\frac{1}{t(1+tC)^2}-\frac{1}{t} \right)=0 .
\eeq 
The term in parenthesis can be written as

\beq \label{eq:secfact}
-\frac{C^2t+2C-\si}{(1+tC)^2} \qquad \mbox{\rm  or } \qquad \frac{(1+t\si)-(1+tC)^2}{t(1+tC)^2}
\eeq
Now, we distinguish two cases. 

Suppose first that for some value of
$t$ the function $1+tC$ vanishes together with its first derivative.
Then 

\beq
M=\frac{1+t\si}{t}\left(\ba{cc} -1 & 1/t \\ t & -1 \ea\right). 
\eeq

As a byproduct, $M_{21}+tM_{22}$ vanishes, and the equivariant Euler
class vanishes. So clearly, the function $1+tC$ cannot vanish
everywhere if we are to find a non-trivial class. Anyway, the
vanishing of $1+tC$ would mean that $e^{\vf}=t_0/t$ for some constant
$t_0$ leading to a metric singular at the origin. It is likely that in
this case, a careful computation with distributions would give a
curvature concentrated at the origin, but we are not interested in
this anyway. 

On the open intervals where $1+tC \neq 0$ the second factor of
(\ref{eq:detm=0}) has to vanish. We get    

\beq
\frac{\si}{(1+tC)^2}+\frac{1}{t(1+tC)^2}-\frac{1}{t}=\frac{1}{t_0}
\eeq
for some constant $t_0$. Using (\ref{eq:secfact}), one obtains

\beq
M=\frac{1+tC}{t_0}\left(\ba{cc} -(1+tC) & (1+(t_0+t)C)(t_0+t)^{-1} \\
(1+tC)(t_0+t)  & -(1+(t_0+t)C) \ea\right)
\eeq
leading to a remarkable simplification of (\ref{eq:eulerclass})~:

\beq \label{eq:final}
E_{int}^{eq}=2^{n/2}\left(\frac{t_0}{t_0+t}\right)^{1/2} \Pf \left(\Om
^{ij}-\frac{1}{t_0+t}\Psi ^i \Psi ^j\right).
\eeq

Now, as $M_{21}=1+t\si$, which has to remain strictly positive, $1+tC$
cannot vanish at the boundary of an open interval where it is nonzero.
This means that $1+Ct$ vanishes nowhere, and that formula
(\ref{eq:final}) is valid everywhere. This is our final formula for the
equivariant Euler class if we decide to trade flexibility (arbitrary
choice of $\vf$ and $\si$) for simplicity ($so(n)$ invariance on $V$,
leading to a simple Pfaffian). The Matha\"{\i}--Quillen representative
never shows up for $n > 2$. 

Some comments are in order. Usually the Thom class is defined by using
function with compact support  (differential topology) or rapid
decrease at infinity (quantum field theory) on $V$. The
Matha\"{\i}--Quillen representative belongs to this second
category. With the general formula, the freedom on $\vf$ and $\si$
allows us to impose any behavior at infinity \footnote{In fact to localize
the Thom class on arbitrary spherical shells if this proves useful.}. On
the other hand our rigid proposal for the Thom class does not decrease
fast at infinity. Despite the fact that this may be
inconvenient in certain applications, we would like to point that it
makes sense nevertheless. To define the Thom class, the crucial point
is that the cohomology of 
$V$ with coefficients having compact support or rapid decrease at infinity is
concentrated in the dimension of $V$ and one dimensional there. It
seems clear that a cohomology of $V$ can be build such as to retain
this property and accept our rigid representative as a well-defined
cohomology class. For instance $k$-forms on $V$
such that for any non-negative integer $l$ the partial derivatives
of order $l$ of the coefficients 
exist and are $O(t^{-\frac{k+l+1}{2}})$ at infinity, endowed with
the usual exterior derivative, should work. 

In particular, we can normalize things in such a way that the integral
on $V$ of the term independent of $\Om $ is 1, as is usual for the
Thom class. A simple calculation gives for the normalized Thom class

\beq
T_V=\frac{1}{(2\pi)^d}\left(\frac{t_0}{t_0+t}\right)^{1/2}
\Pf \left(\frac{1}{t_0+t}\Psi ^i \Psi ^j - \Om ^{ij}\right).
\eeq 
Playing with the value of $t_0$ allows to localize around the zero section.
This formula already appears in \cite{hl93} as a specialization of another
formula for the Thom class.

\section{Conclusion.}

In these notes, we have obtained formul\ae\ for the universal Thom
class of a vector bundle. A special choice leads to a rigid
representative involving Cauchy-type kernels. 
It would be very interesting to know whether the Matha\"{\i}--Quillen
representative, with its Gaussian-type kernel, is also a rigid member
of some natural family of representatives.

\noindent

\section*{\large \bf Acknowledgments}

\noindent We are very grateful to Raymond Stora and Georges Girardi
for many fruitful discussions. We thank Reese Harvey for making us aware of
ref. \cite{hl93}, which contains material closely related to this work.


\begin{thebibliography}{BCI94}

\bibitem[Bl92]{Bl92} M. Blau, ''The Matha\"{\i} Quillen formalism and
topological field theory'', Lectures at the 28th Karpacz Winter School on
Infinite dimensional geometry in physics, Karpacz, Poland, Feb. 1992,
NIKHEF-H/92-07, 17-29

\bibitem[BGV91]{BGV91} N. Berline, E. Getzler and M. Vergne, ''Heat
Kernels and Dirac Operators'', Grundlehren des Mathematischen
Wissenschasft 298, Springer-Verlag Berlin Heidelberg 1992

\bibitem[BS88]{BS88} L. Baulieu and I.M. Singer, Nucl. Phys. B Proc Supl
15B (1988) 12

\bibitem[BS91]{BS91} L. Baulieu and I.M. Singer, Commun. Math. Phys. 135
(1991) 253

\bibitem[BT82]{BT82} R. Bott and L.W. Tu, Differential Forms in Algebraic
Topology , Springer-Verlag 1982

\bibitem[C50]{C50} H. Cartan, Colloque de Topologie (Espaces Fibr\'{e}s),
Brussels 1950, CBRM, 15-56

\bibitem[DV93]{DV93} M. Duflo and M. Vergne, Sur la cohomologie \'{e}quivariante des vari\'{e}t\'{e}s diff\'{e}rentiables, Ast\'{e}risque 215, 1993, 5-108

\bibitem[HL93]{hl93} R. Harvey and B. Lawson, A theory of characteristic
currents associated with a singular connection", Ast\'erisque 213 (1993),
1-268.

\bibitem[K93]{K93} J. Kalkman, Commun. Math. Phys. 153 (1993) 447

\bibitem[MQ86]{MQ86} V. Matha\"{\i} and D. Quillen, Topology 25 (1986) 85

\bibitem[MSZ85]{MSZ85} J. Manes, R. Stora and B. Zumino, Commun. Math.
Phys. 102 (1985) 157

\bibitem[OSB89]{OSB89} S. Ouvry, R Stora and P. van Baal, Phys. Lett.
B220
(1989) 590

\bibitem[STW94]{STW94} R Stora, F. Thuillier and J.C. Wallet, Lectures at
the 1st Caribbean Spring School of Mathematics and Theoretical Physics,
Saint Fran\c {c}ois, Guadeloupe, May 30 - June 13 1993, Proceedings 1995
(preprint ENSLAPP-A-481/94).

\bibitem[W88]{W88} E. Witten, Commun. Math. Phys. 117 (1988) 353

\bibitem[WBS88]{WBS88} E. Witten, Commun. Math. Phys. 118 (1988) 411. L.
Baulieu and I.M. Singer, Commun. Math. Phys. 125 (1989) 227
\end{thebibliography}
\end{document}